\documentclass[a4paper,10pt]{amsart}

\usepackage{amsmath,amssymb,amsthm,a4wide}

\usepackage{tikz}
\usetikzlibrary{shapes}
\usetikzlibrary{intersections}

\usepackage{graphicx}
\usepackage{graphics}

\newtheorem*{thm}{Theorem}

\newtheorem*{corollary}{Corollary}
\newtheorem*{proposition}{Proposition}
\newcommand{\var}{\operatorname{var}}

\begin{document}

\title[]{An uncertainty principle on compact manifolds}
\author{Stefan Steinerberger}

\address{Department of Mathematics, Yale University, 10 Hillhouse Avenue, New Haven, CT 06511, USA}
\email{stefan.steinerberger@yale.edu}

\begin{abstract}Breitenberger's uncertainty principle on the torus $\mathbb{T}$ and its higher-dimensional analogue on $\mathbb{S}^{d-1}$ are well understood. We give describe an entire family of uncertainty principles on compact manifolds $(M,g)$, which includes the classical Heisenberg-Weyl uncertainty principle (for $M=B(0,1) \subset \mathbb{R}^d$ the unit ball with the flat metric) and the Goh-Goodman uncertainty principle (for $M=\mathbb{S}^{d-1}$ with the canonical metric) as special cases. This raises a new geometric problem related to small-curvature low-distortion embeddings: given a function $f:M \rightarrow \mathbb{R}$, which uncertainty principle in our family yields the best result? We give a (far from optimal) answer for the torus, discuss disconnected manifolds and state a variety of other open problems.
\end{abstract}
\maketitle

\section{Introduction}
\subsection{Classical uncertainty principle.} The classical uncertainty principle on $\mathbb{R}$ needs no introduction: for $L^2-$normalized Schwartz functions we may use partial integration and get
$$ 1 = \int_{\mathbb{R}}{u^2dx} = \int_{\mathbb{R}}{\left(\frac{d}{dx}x\right)u^2dx} 
= -\int_{\mathbb{R}}{x 2u u_x dx}  \leq 2 \left(\int_{\mathbb{R}}{u^2 x^2 dx}\right)^{\frac{1}{2}}
\left(\int_{\mathbb{R}}{u_x^2dx}\right)^{\frac{1}{2}}$$
with equality for Gaussians; generalizations to $\mathbb{R}^d$ are immediate. Indeed, the way this
uncertainty principle is phrased here (with a spatial derivative instead of the Fourier transform),
it is immediately clear that certains versions of it will certainly be true on suitable manifolds: one can replace $|x|$ by the geodesic distance $d(x, x_0)$ to a fixed point $x_0 \in M$ and all other 
quantities are still well-defined on a manifold. In the non-compact case, the growth of geodesic balls (and thus curvature) of the manifold will play an important role; a typical question in this line of research is not so much concerned with the validity of an uncertainty principle but finding the sharp constant, see for example \cite{ko1, ko2} or the survey of Folland \& Sitaram \cite{foll}. 

\subsection{Breitenberger's uncertainty principle.} Breitenberger was inspired by the idea that 'uncertainty
measures must no depend on the choice of the origin of the measurement scale' \cite{breit} and, in the
case of manifolds, should not depend on a distinguished point $x_0 \in M$ -- a simple form of such 
an inequality would be
$$ \inf_{p \in M}\left(\int_{M}{u^2 d(x,p)^2 dg}\right)\left(\int_{M}{|\nabla u|^2dg}\right) \gtrsim \|u\|_{L^2(M)}^4.$$
Breitenberger \cite{breit} pointed out that on the torus $\mathbb{T}$ there exists a natural inequality that implies even more: if 
$$f = \sum_{k = -\infty}^{\infty}{c_k e^{i k x}} \qquad \mbox{and} \qquad
\|f\|_{L^2(\mathbb{T})} = 1,$$
then we may define the mean localization $\tau(f) \in \mathbb{C}$ as
$$ \tau(f) = \frac{1}{2\pi} \int_{\mathbb{T}}{e^{i x}|f(x)|^2 dx} \in \mathbb{C}.$$
If we furthermore define the frequency variance $\var_F(f) \in \mathbb{R}$ as
$$ \var_{F}(f) = \left(\sum_{k = -\infty}^{\infty}{k^2|c_k|^2}\right)
- \left( \sum_{k = -\infty}^{\infty}{k|c_k|^2}\right)^2$$
and the angular variance $\var_A(f)$ as
$$ \var_{A}(f)  = \frac{1-|\tau(f)|^2}{|\tau(f)|^2},$$
then we may state Breitenberger's uncertainty principle as
$$ \var_F(f)\var_A(f) > \frac{1}{4},$$
whenever all quantities are defined. The uncertainty principle has two different regimes of interest: if $0 <|\tau(f)| \ll 1$, then $\var_A(f) \sim 1$ and the inequality provides a lower bound on $\var_F(f)$ and tells us that the gradient cannot be too small: if it was very small, the function would have to be almost constant and $\tau(f)$ would be much closer to $0$. The second statement is perhaps more interesting: if $|\tau(f)|$ is very close to 
1, then $0 \leq \var_A(f) \ll 1$ and the function $f$ has to have its $L^2-$mass concentrated
around a point, which requires it to have a large gradient.
\begin{center}
\begin{figure}[h!]
\begin{tikzpicture}[scale = 1.7]
\draw [ultra thick] (0,0) circle [radius=1];
\draw [ultra thick, fill] (0.1,-0.1) circle [radius=0.03];
\draw (0.35,-0.2) node{\Large $\tau(f)$};
\draw [->] (0,0) -- (1.2,0);
\draw [->] (0,0) -- (0,1.2);
\draw (-1.2,0) -- (0,0);
\draw (0,-1.1) -- (0,0);
\draw (1.5,-1.3) node{\Large $\tau(f) = \int_{\mathbb{T}}{(\cos(t),\sin(t))|f(t)|^2dt}$};
\draw [ultra thick] (3,0) circle [radius=1];
\draw [ultra thick, fill] (3.6,-0.5) circle [radius=0.03];
\draw (3.3,-0.65) node{\Large $\tau(f)$};
\draw [->] (3,0) -- (4.2,0);
\draw [->] (3,0) -- (3,1.2);
\draw (3-1.2,0) -- (3,0);
\draw (3,-1.1) -- (3,0);
\end{tikzpicture}
\caption{Two regimes: (left) if the center of mass is not in the origin, then the function $f$ is not constant and (right) if the center of mass is close to the boundary, the function $f$ is tightly concentrated.}
\end{figure}
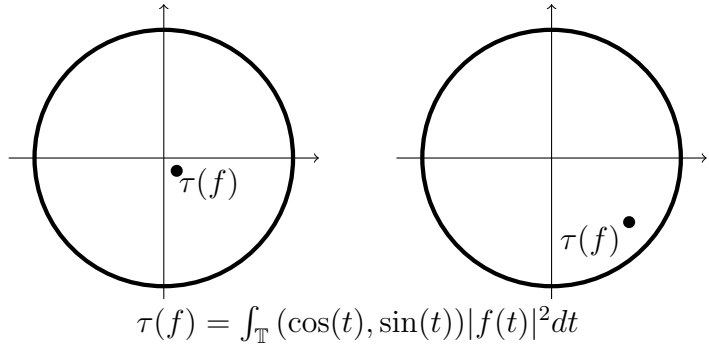
\end{center}
\vspace{-15pt}
The two regimes displayed in the picture are indeed fundamentally different: roughly speaking, the
regime $0 < |\tau(f)| \ll 1$ behaves like a Poincar\'{e}-type inequality while $|\tau(f)| \sim 1$ behaves
more like the classical Euclidean uncertainty principle. 

\subsection{Further work.} Prestin \& Quak \cite{pr} have shown that the constant
in the Breitenberger uncertainty principle is sharp. Moreover, Prestin, Quak, Rauhut \& Selig \cite{all4} have proven that if we take a function
defined on the real line and periodize it in a suitable way, then the Breitenberger uncertainty principle reduces to the classical Euclidean uncertainty principle in the limit. There exists a natural generalization on the sphere $\mathbb{S}^d$ which has been of considerable interest.
\begin{thm}[Goh \& Goodman \cite{goh}] Let $g$ be the normalized surface measure on $\mathbb{S}^d$,
let $f \in H^1(\mathbb{S}^d, g)$ be $L^2-$normalized. If we define
$$ \tau(f) = \int_{\mathbb{S}^d}{x|f(x)|^2 dg(x)},\quad  \var_{A}(f)  = \frac{1-|\tau(f)|^2}{|\tau(f)|^2} \quad \mbox{and} \quad \var_F(f) = \left\langle -\Delta_{\mathbb{S}^d} f, f\right\rangle,$$
then
$$ \var_F(f)\var_A(f) \geq \frac{d^2}{4}.$$
\end{thm}
This result generalizes earlier theorems of Narcowich \& Ward \cite{nar} as well as R\"osler \& Voit \cite{ros} and is now considered the canonical version of an uncertainty principle on the sphere. 
It has one considerable disadvantage: whenever $|f(x)|$ is symmetric around one of the coordinate axes, then $\tau(f) = \textbf{0} \in \mathbb{R}^{d+1}$ and $\var_{A}(f)$ is undefined and no
lower bound on $\var_{F}(f)$ can be recovered. This is not at all surprising since constant functions have the same symmetry and their gradient vanishes but raises the question of whether this
obstruction could somehow be avoided.\\

The only general result is due to Erb \cite{erb}, who considers arbitrary compact, connected manifolds $(M,g)$ without boundary, a distinguished point $p \in M$ and functions radial around that point. This gives a highly nontrivial generalization of the Breitenberger uncertainty principle for a large class of differential operators and weights; however, due to the condition of radiality, the result can also be regarded as a weighted one-dimensional result. We were motivated by the lack of a genuine $n-$dimensional version of the inequality and the above mentioned difficulty arising whenever functions exhibit a symmetry.\\

\subsection{Our approach.} The concept of frequency variance $\var_{F}(f)$ series is intimately tied to the trigonometric expansion on $\mathbb{T}$.  As was also done in the above inequality for the sphere,
we will proceed as in the the classical Euclidean case and use
$$ \int_{M}{|\nabla f|^2 dg} \qquad \mbox{as a measure of frequency localization.}$$
Identifying $\mathbb{C} \cong \mathbb{R}^2$ and
using $m: \mathbb{T} \rightarrow \mathbb{R}^2$ to denote the map $$m(x) = (\cos(x), \sin(x)),$$ we may rewrite the mean localization in a more geometric way
$$ |\tau(f)|^2 = \left|\frac{1}{2\pi}\int_{\mathbb{T}}{e^{i x}|f(x)|^2 dx}\right|^2 = \left\| 
\frac{1}{2\pi}\int_{\mathbb{T}}{m(x)|f(x)|^2 dx}\right\|^2_{\ell^2(\mathbb{R}^2)},$$
where $\ell^2(\mathbb{R}^2)$ denotes the Euclidean norm in the plane. Replacing the frequency variance by
the gradient, Breitenberger's uncertainty principle implies the already mentioned inequality on the sphere for $d=1$ (where we identify $\mathbb{T} \cong \mathbb{S}^1$)
$$ \frac{1-|\tau(f)|^2}{|\tau(f)|^2} \int_{\mathbb{T}}{|\nabla f|^2dg} \geq \frac{1}{4} \quad
\mbox{for all}~f \in H^1(\mathbb{T})~\mbox{with}~\|f\|_{L^2(\mathbb{T})} = 1.$$
The advantage of this formulation is that all these terms can be interpreted in a meaningful way
for arbitrary compact manifolds embedded in some Euclidean space. Our next step is slight reformulation of the statement, where we replace $1-|\tau(f)|^2$ in the numerator by $1-|\tau(f)|$. This has no major impact on the statement itself and only changes the constant:  whenever $|\tau(f)|$ is very small, the precise form of that term is not important as it is dominated by the constant 1 anyway and whenever $|\tau(f)| = 1-\varepsilon$, for $\varepsilon$ small, then $$ 1-|\tau(f)|^2 = 2\varepsilon - \varepsilon^2 \sim \varepsilon = 1-|\tau(f)|.$$
Our motivation for considering this modified inequality
$$ \frac{1-|\tau(f)|}{|\tau(f)|^2} \int_{\mathbb{T}}{|\nabla f|^2dg} \geq \frac{1}{8} \quad
\mbox{for all}~f \in H^1(\mathbb{T})~\mbox{with}~\|f\|_{L^2(\mathbb{T})} = 1$$
is that it will be the template for our main theorem: our general uncertainty principle will reduce 
to that very inequality in the case of $(M,g) = (\mathbb{T},~$can$)$ and $m:\mathbb{T} \rightarrow \mathbb{R}^2$ given by $m(x) = (\cos{x}, \sin{x})$. The change of the numerator 
is motivated by the fact that our result should have the same scaling as the classical uncertainty
principle in Euclidean space. 

\subsection{Setup} Let $(M,g)$ be some compact, connected $n-$dimensional, smooth manifold (without boundary or with smooth boundary). We will study the image of the manifold $M$ under some suitable embedding into Euclidean space $m: M \rightarrow \mathbb{R}^d$. Each admissible mapping $m$ will give
a particular uncertainty principle on $(M,g)$. This allows a greater flexibility and will allow us to avoid the aforementioned problems arising when functions exhibit a symmetry.
We demand that the maps $m: M \rightarrow \mathbb{R}^d$ satisfy three properties:
\begin{enumerate}
\item  \textit{Bilipschitz.}
 $m$ is bilipschitz: there is a constant $L < \infty$ for which
$$\forall~x,y \in M: \quad  \frac{d_g(x,y)}{L} \leq \|m(x) - m(y)\|_{\ell^2(\mathbb{R}^d)} \leq  Ld_g(x,y),$$
where $d_g(\cdot, \cdot)$ is the geodesic distance on the manifold.

\item \textit{Curvature}. $m$ satisfies the following curvature condition for some $0 < C < \infty$, for all
integers $N \geq 1$ and all elements $x_1, \dots, x_N,z \in M$:
$$ \left\| \frac{1}{N}\sum_{i=1}^{N}{m(x_i)} - m(z)\right\|_{\ell^2(\mathbb{R}^d)} \geq 
 \frac{C}{N}\sum_{i=1}^{N}\left\|m(x_i) - m(z) \right\|^2_{\ell^2(\mathbb{R}^d)}.$$
\item \textit{Normalization.} The third condition on $m$ can always be satisfied w.l.o.g. and merely eliminates the translation invariance of the problem to simplify notation: henceforth, we demand that
$$ \int_{M}{m(x) dg} = \textbf{0} \in \mathbb{R}^d.$$
\end{enumerate}
The second condition implicitely demands some curvature by saying that no linear combination of elements
in $ \left\{m(x): x \in M\right\} \subset \mathbb{R}^d$ lies close to another element of $ \left\{m(x): x \in M\right\}$ unless all elements were chosen very close to that element to begin with.
\begin{figure}[h!]
\begin{tikzpicture}[scale = 1.3]
\draw[very thick] (0,0) to [out=00,in=240] (3,1);
\draw [fill] (0,0) circle [radius=0.05];
\draw [fill] (3,1) circle [radius=0.05];
\draw [fill] (1.5,0.5) circle [radius=0.05];
\draw [fill] (1.63,0.08) circle [radius=0.05];
\draw (1.8,-0.1) node{$m(z)$};
\draw [dashed] (1.5,0.5) -- (1.63,0.08);
\draw (0,0) -- (3,1);
\draw (-0.3,-0.3) node{$m(x)$};
\draw (3.4,0.9) node{$m(y)$};
\draw (1.2,0.8) node{\Large $\frac{m(x)+m(y)}{2}$};
\end{tikzpicture}
\caption{The geometric condition illustrated for $N=2$.}
\end{figure}
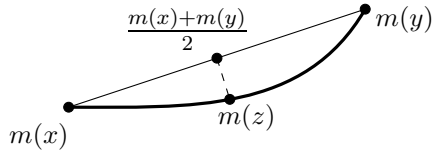
 The condition could be explicitely phrased in terms of curvature in simple special cases: for example, if we are looking at admissible mappings $m:\mathbb{T} \rightarrow \mathbb{R}^2$,
then the second condition requires that the curvature of $m(\mathbb{T})$ is uniformly bounded away from 0. The way the condition is formulated, however, does allow for the curvature to be
unbounded (i.e. $m(\mathbb{T})$ need not be everywhere differentiable). We emphasize that the second condition is a statement about the embedding $ \left\{m(x): x \in M\right\}$
and not any notion of curvature defined on the manifold $M$ itself.
For the type of estimates we seek, all three conditions are necessary.
\begin{center}
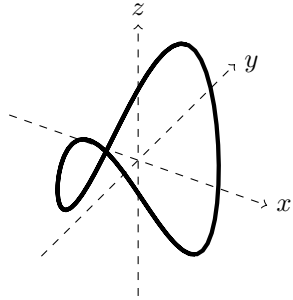
\begin{figure}[h!]
\begin{tikzpicture}[scale= 0.9, x={(0.707cm,0.707cm)},z={(0cm,1cm)},y={(-0.942cm,0.3333cm)}]
\draw[dashed, ->] (-2,0,0) -- (2,0,0) node[right] {$y$};
\draw[dashed, ->] (0,2,0) -- (0,-2,0) node[right] {$x$};
\draw[dashed, ->] (0,0,-2) -- (0,0,2) node[above] {$z$};
\draw [ultra thick] (1,0,1)
\foreach \z in {0,0.1,...,10}
{ -- ({cos(\z*2*3.1415*10)},{sin(\z*2*3.141*10)},{cos(\z*2*2*3.1415*10)})};
\end{tikzpicture}
\caption{An admissible embedding $m:\mathbb{T} \rightarrow \mathbb{R}^3$.}
\end{figure}
\vspace{-20pt}
\end{center}
The simplest example of such a mapping is $(M,g) = (\mathbb{S}^d, $~can$)$ and $m$ being the canonical embedding. Other geometric examples are given by the boundary of the unit ball of $\ell_p^n$ for any $1 < p < \infty$ or the graph of any strictly convex function defined on a unit ball $B(0,1) \subset \mathbb{R}^d$. A slightly more
nontrivial example is given by $(M,g) = (\mathbb{T}, $~can$)$ and $m:\mathbb{T} \rightarrow \mathbb{R}^3$ given by
$$ m(t) = (\cos{t}, \sin{t}, \phi(t))$$
for any smooth $\phi:\mathbb{T} \rightarrow \mathbb{R}$. It is easy to see by considering the projection
onto the first two components that the statement holds for some $C_{\phi}$ depending on $\phi$. 
Another nontrivial example is given by the two-dimensional manifold $M = \mathbb{T} \times [0,1]$ for which one needs a properly chosen map $m:M \rightarrow \mathbb{R}^3$ so that the properties are satisfied.
\begin{center}
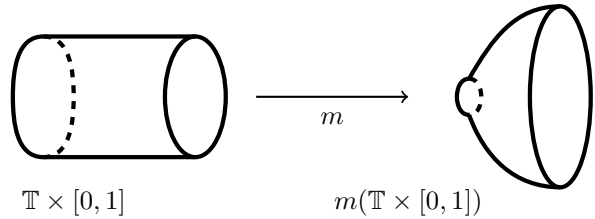
\begin{figure}[h!]
\begin{tikzpicture}[scale = 0.4]
\draw [ultra thick] (-2,0) ellipse (1cm and 2cm);
\draw [ultra thick] (-7,2) to [out=180,in=90] (-8,0)  to [out=270,in=180] (-7,-2);
\draw [ultra thick, dashed] (-7,2) to [out=0,in=90] (-6,0)  to [out=270,in=0] (-7,-2);
\draw [ultra thick] (-7,2) -- (-2,2);
\draw [ultra thick] (-7,-2) -- (-2,-2);
\draw (-6,-3.5) node{$\mathbb{T} \times [0,1]$};

\draw [thick] [->] (0,0) -- (5,0);
\draw (2.5,-0.7) node{$m$};
\draw [ultra thick] (10,0) ellipse (1cm and 3cm);
\draw [ultra thick] (7,0.6) to [out=180,in=90] (6.6,0)  to [out=270,in=180] (7,-0.6);
\draw [ultra thick, dashed] (7,0.6) to [out=0,in=90] (7.4,0)  to [out=270,in=0] (7,-0.6);
\draw [ultra thick] (10,3) to [out=180,in=60] (7,0.6);
\draw [ultra thick] (10,-3) to [out=180,in=290] (7,-0.6);
\draw (5,-3.5) node{$m(\mathbb{T} \times [0,1])$};
\end{tikzpicture}
\caption{Blowing up $\mathbb{T} \times [0,1]$ in a suitable way.}
\end{figure}
\vspace{-15pt}
\end{center}

\subsection{Statement.} We now give our main statement. We will use the notation $A \lesssim_{x} B$ to denote the existence of an absolute constant depending only on $x$. If both $A \lesssim_{x} B$ and $A \gtrsim_{x} B$, then we write $A \sim_{x} B$. Whenever a line contains the letters $C, L$, these letters are too be understood as the constant of the mapping $m$ occuring
in the same line (which is simpler than writing $C(m), L(m)$). \\

Our main statement is a Breitenberger-type uncertainty principle for an arbitrary map $m:\mathcal{M} \rightarrow \mathbb{R}^d$. In particular, if we restrict to maps which are not too degenerate, i.e.
$$ \left\{m:M \rightarrow \mathbb{R}^d: m~\mbox{admissible and}~ C(m) \leq C \wedge L(m) \leq L\right\},$$
then our inequality is \textit{uniform} over that entire set. 
\begin{thm} Let $(M,g)$ and the map $m: M \rightarrow \mathbb{R}^d$ be admissible. For real-valued $f \in H^1(M)$ with $\|f\|_{L^2(M)} = 1$, we have
$$ \left(\inf_{z \in M}{\left\|m(z) - \int_{M}{m(x)|f(x)|^2dg}\right\|_{\ell^2(\mathbb{R}^d)}}\right) \left\| \int_{M}{m(x)|f(x)|^2dg}\right\|_{\ell^2(\mathbb{R}^d)}^{-2} \int_{M}{|\nabla f|^2dg} \gtrsim_{(M,g),d} \frac{C}{L^4}.$$
We emphasize again that the lower bound does not depend on the function $f$ nor on any other properties of $m:M \rightarrow \mathbb{R}^d$ except the size of $C,L$ and $d$. 
\end{thm}
This should be compared to the uncertainty principle on the sphere $\mathbb{S}^{d-1}$ or the uncertainty principle on $\mathbb{T}$: then the first term is merely $1-|\tau(f)|$ while the second
term is $|\tau(f)|^{-2}$. Both terms continue playing the same role but the lack of spherical symmetry of $m(M)$ requires us to write the first term in a more general way.\\

The right-hand side has the sharp
scaling in $C$ and $L$ (see below for an example). 
 The statement has a straightforward extension to the case of disconnected manifolds; we are going to
establish that result in Section 4.
 If $(M,g) = (B(0,1),~$flat$)$
is the flat unit ball in $\mathbb{R}^d$, then the main statement implies, for a suitable choice of $m$, the classical uncertainty principle in $\mathbb{R}^d$ (this will be described in greater detail in Section 2).
For $(M,g) = (\mathbb{S}^d, $~can$)$ and $m$ being the canonical embedding the statement simplifies as follows.
\begin{corollary} Let $f\in H^1(\mathbb{S}^{d-1})$ be real-valued with the normalization $\|f\|_{L^2(\mathbb{S}^{d-1})} = 1$. Then
$$ \left(1-\left\|\int_{\mathbb{S}^{d-1}}{x|f(x)|^2dg(x)} \right\|_{\ell^2(\mathbb{R}^{d})}\right)\left\|\int_{\mathbb{S}^{d-1}}{x|f(x)|^2dg(x)} \right\|_{\ell^2(\mathbb{R}^{d})}^{-2}\int_{\mathbb{S}^{d-1}}{|\nabla f|^2dg} \gtrsim_d 1.$$
\end{corollary}
This is almost precisely (up to a power in the first term that only affects the constant) the inequality studied by Narcowich \& Ward \cite{nar} (for $d=3$), R\"osler \& Voit \cite{ros} (for arbitrary $d$ and radial functions) and, unconditionally, Goh \& Goodman \cite{goh}. The work of Goh \& Goodman implies that the optimal
constant $c_d$ in the inequality satisfies $(d-1)^2/8 \leq c_d \leq (d-1)^2/4.$ However, as already emphasized above, our result is much more general than that
and one could also consider mappings $m:M \rightarrow \mbox{(some ellipsoid)}$ or another shape entirely; a special case that is easy to state is that of
$\ell^p_n-$balls. It follows from Clarkson's inequalities \cite{clark} that the unit ball in $\ell^p_n$
$$ B_{p,d} =  \left\{x \in \mathbb{R}^d: \|x\|_{p} = 1\right\} \qquad \mbox{is uniformly convex if}~p>1.$$
By taking a suitable map $m:\mathbb{S}^{d-1} \rightarrow B_{p,d}$, we get the following immediate generalization of the above corollary.
\begin{corollary} Let $f\in H^1(\mathbb{S}^{d-1})$ be real-valued with the normalization $\|f\|_{L^2(\mathbb{S}^{d-1})} = 1$. Then, for every $p>1$, 
$$ \left(1-\left\|\int_{\mathbb{S}^{d-1}}{m(x)|f(x)|^2dg(x)} \right\|_{\ell^p(\mathbb{R}^{d})}\right)\left\|\int_{\mathbb{S}^{d-1}}{m(x)|f(x)|^2dg(x)} \right\|_{\ell^2(\mathbb{R}^{d})}^{-2}\int_{\mathbb{S}^{d-1}}{|\nabla f|^2dg} \gtrsim_{p,d} 1.$$
\end{corollary}
Using the equivalence of norms in finite-dimensional vector spaces, one could additionally replace $\ell^2(\mathbb{R}^d$ by $\ell^q(\mathbb{R}^d)$ for any $1 \leq q \leq \infty$. We note that the restriction
$p > 1$ is necessary because the $\ell^1_n$ contains flat segments, which can be used to construct a counterexample.

\subsection{Sharp scaling of the constant.} The purpose of this short section is to demonstrate that the scaling $C/L^4$ obtained for the right-hand side of the main statement is optimal. We consider a particular example $m:\mathbb{S}^1 \rightarrow 
\mathbb{R}^2$, partition the unit sphere $\mathbb{S}^1$ into four parts of equal length, fix some large constant $L \gg 1$ and map it to a circle with radius $\sim L$ by 
shrinking two intervals while expanding two others. We demand that the map is symmetric and that in particular
\begin{align*}
 \|m(a)-m(b)\|_{\ell^2(\mathbb{R}^2)} &= \|m(a)-m(b)\|_{\ell^2(\mathbb{R}^2)}  \sim L^{-1} \\
 \|m(a)-m(c)\|_{\ell^2(\mathbb{R}^2)} &= \|m(b)-m(d)\|_{\ell^2(\mathbb{R}^2)}  \sim L.
\end{align*} 
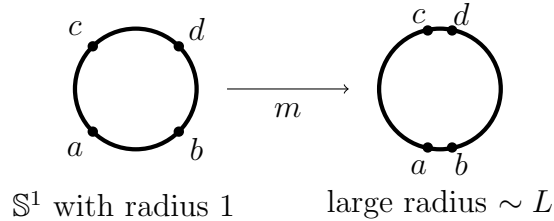
\begin{center}
\begin{figure}[h!]
\begin{tikzpicture}[scale = 0.8]
\draw [ultra thick] (0,0) circle [radius=1];
\draw [ultra thick, fill] (-0.707,-0.707) circle [radius=0.05];
\draw [ultra thick, fill] (0.707,-0.707) circle [radius=0.05];
\draw [ultra thick, fill] (0.707,0.707) circle [radius=0.05];
\draw [ultra thick, fill] (-0.707,0.707) circle [radius=0.05];
\draw (-1,-1) node{\Large $a$};
\draw (1,-1) node{\Large $b$};
\draw (-1,1) node{\Large $c$};
\draw (1,1) node{\Large $d$};
\draw (-0.2,-1.9) node{\Large $\mathbb{S}^1$ with radius 1};
\draw [->] (1.5,0) -- (3.5,0);
\draw (2.5,-0.35) node{\Large $m$};

\draw [ultra thick] (5,0) circle [radius=1];
\draw [ultra thick, fill] (5-0.2,-0.97) circle [radius=0.05];
\draw [ultra thick, fill] (5.2,-0.97) circle [radius=0.05];
\draw [ultra thick, fill] (5.2,.97) circle [radius=0.05];
\draw [ultra thick, fill] (5-0.2,.97) circle [radius=0.05];
\draw (5-0.35,-1.3) node{\Large $a$};
\draw (5+0.35,-1.25) node{\Large $b$};
\draw (5-0.35,1.2) node{\Large $c$};
\draw (5+0.35,1.25) node{\Large $d$};
\draw (5,-1.9) node{\Large large radius $\sim L$};
\end{tikzpicture}
\end{figure}
\end{center}
Since the image of $m$ is merely a circle with radius $\sim L$, the constant $C$ behaves
like the infimum of the curvature and satisfies $C \sim L^{-1}$. 
Take now some $f:\mathbb{S}^1 \rightarrow \mathbb{R}$,
which is supported in the interval between $a$ and $b$, vanishes outside and satisfies
$$ \|f\|_{L^2(\mathbb{T})} = 1 \quad \mbox{and} \quad \|\nabla f\|_{L^2(\mathbb{T})} \sim 1.$$
A simple calculation shows that
\begin{align*}
\inf_{z \in M}{\left\|m(z) - \int_{M}{m(x)|f(x)|^2dg}\right\|_{\ell^2(\mathbb{R}^d)}} &\sim \frac{1}{L^3} \\
\left\|\int_{M}{m(x)|f(x)|^2dg}\right\|^{-2}_{\ell^2(\mathbb{R}^d)} &\sim \frac{1}{L^2}
\end{align*}
and thus the uncertainty principle scales as $L^{-5}$ matching the right hand side $C/L^4 \sim L^{-5}$.
For any desired $C \ll L^{-1}$ one can additionally construct another example realizing $\sim C/L^4$ 
precisely by replacing the arcs between $a,b$ and $c,d$ in a symmetric fashion by a suitable paraboloid.

\subsection{Uniform families and an inverse problem.} Given now $(M,g)$ and a function $f:M \rightarrow \mathbb{R}$, what is the
mapping $m:M \rightarrow \mathbb{R}^d$ which allows to recover the best lower bound on the derivative? Recall that our main statement
reads as follows.
$$ \left(\inf_{z \in M}{\left\|m(z) - \int_{M}{m(x)|f(x)|^2dg}\right\|_{\ell^2(\mathbb{R}^d)}}\right) \left\| \int_{M}{m(x)|f(x)|^2dg}\right\|_{\ell^2(\mathbb{R}^d)}^{-2} \int_{M}{|\nabla f|^2dg} \gtrsim_{(M,g),d} \frac{C}{L^4}.$$
Using the trivial estimate,
$$ \inf_{z \in M}{\left\|m(z) - \int_{M}{m(x)|f(x)|^2dg}\right\|_{\ell^2(\mathbb{R}^d)}} \lesssim_{(M,g)} L,$$
we may simplify (and weaken) our inequality in a way that puts a bigger emphasis on the actual inverse problem
$$  \int_{M}{|\nabla f|^2dg} \gtrsim_{(M,g),d} \frac{C}{L^5}\left\| \int_{M}{m(x)|f(x)|^2dg}\right\|_{\ell^2(\mathbb{R}^d)}^{2}.$$
We are interested in how one would actually chose $m$ to get the best lower bound on $\|\nabla f\|_{L^2}^2$. Geometrically, it corresponds
to a building an embedding that is suitably taylored to $f$ -- it is not clear to us how much of $\|\nabla f\|_{L^2}^2$ can actually be recovered
or how one would find such maps.
\begin{quote} \textbf{Inverse problem.} Given $(M,g)$ and $f:M \rightarrow \mathbb{R}$, which map $m:M \rightarrow \mathbb{R}^d$ is both very curved ($C$ big) as well as
isometric as possible ($L$ small) but also has the property that when it is weighted with $|f(x)|^2$ yields a point far away from the origin? In short,
prove lower bounds
$$\sup_{m:M \rightarrow \mathbb{R}^d}\frac{C}{L^5}\left\| \int_{M}{m(x)|f(x)|^2dg}\right\|_{\ell^2(\mathbb{R}^d)}^{2} \qquad \mbox{in terms of}~f$$
and describe the mappings $m$ achieving that bound.
\end{quote}

The following is obvious: any such lower bound needs to decay as $\|\nabla f\|_{L^{\infty}}$ increases since a large
gradient allows for an oscillating function whose $L^2-$norm is evenly distributed over a
manifold: a mapping $m$ with not too large distortion $L$ does not 'see' the oscillations and treats such
a function like it would the constant function. On the other hand, any reasonable bound needs to incorporate
the fact that the quantity vanishes if $f$ is constant.\\

 The most natural case to study is again $m:\mathbb{T} \rightarrow \mathbb{R}^2$ which already seems very nontrivial. Since there are more maps $m:\mathbb{T} \rightarrow \mathbb{R}^3$, it becomes easier to give \textit{some} nontrivial bound in that case: we note one particularly simple construction. It displays all the characteristics described above but is probably far from optimal.
\begin{proposition} Let $f \in H^1(\mathbb{T})$ be normalized in $L^2(\mathbb{T})$. We have
$$\sup_{m:\mathbb{T} \rightarrow \mathbb{R}^3}\frac{C}{L^5}\left\|\int_{\mathbb{T}}{m(x)|f(x)|^2dg} \right\|_{\ell^2(\mathbb{R}^{3})} \gtrsim \frac{1}{(1+\|\nabla f\|_{L^{\infty}})^7}\left( \int_{\mathbb{T}}{|f(x)|^4dx}
- \frac{1}{2\pi}\right).$$
\end{proposition}
Note that H\"older's inequality implies that the right-hand side never vanishes unless $f$ is constant. We have no reason to believe that this is optimal and
lack any approach to the general problem.

\subsection{Outline.} The remainder of the paper is structured as follows: in the next section we describe an interesting special case and how it implies the classical Euclidean uncertainty principle. Section 3 gives a proof of the main statement, Section 4 discusses an extension of the result to disconnected manifolds, Section 5 describes the arising inverse problem one gets when trying to find the \textit{optimal} map $m$ in a given situation well as some further results and open problems.

\section{The Euclidean uncertainty principle}
The purpose of this section is first to show how the classical Euclidean uncertainty principle follows rather immediately as a special case of our inequality (where we use the symmetries of the Euclidean uncertainty principle to compactify $\mathbb{R}^d$); the second point we will be making is that a more natural statement is also
implicitely encoded in our inequality. We hope that this serves as a nice example clarifying the precise
role of each of the three terms and their interplay and underlines the connection between
our result and Breitenberger's original intuition.

\subsection{Classical uncertainty.} We want to demonstrate that the inequality
$$\left(\int_{\mathbb{R}^d}{|x|^2 f^2 dx}\right)\left(\int_{\mathbb{R}^d}{|\nabla f|^2dx}\right) \gtrsim_d \|f\|_{L^2(\mathbb{R}^d)}^4$$
is implicitely contained in our main theorem; this will not give a new proof of the Euclidean uncertainty principle since we will actually
be using it in the proof of our main theorem, however, it is nonetheless a nice fact that the
main result still contains the most basic uncertainty principle.
\begin{proof}We may use the invariance of the inequality under dilation with some $\lambda > 0$
$$ f(x) \rightarrow \lambda^{\frac{d}{2}}f(\lambda x)$$
 to assume that the support of $f$ is contained in the unit disk $B(0,1) \subset \mathbb{R}^d$, has 99\% of its $L^2-$mass at distance at most $1/10$ from the origin and vanishes on the boundary (this step may be regarded as a compactification of the Euclidean space). We may further assume that $\|f\|_{L^2(B(0,1))} = 1.$ The Polya-Szeg\H{o} inequality on symmetric decreasing rearrangement (see, for example, Lieb \& Loss \cite{lieb}) implies that it suffices to consider radial functions. It remains to show that the gradient term increases
fast enough when a function concentrates its $L^2-$mass around the origin. We choose the manifold to be the unit ball $M = B(0,1) \subset \mathbb{R}^d$ and $g$ to be the flat metric. For reasons that will soon be apparent, we set
$m:M \rightarrow \mathbb{R}^{d+1}$ given by
$$ m(x) = \left(x,|x|^2- \frac{1}{|B(0,1)|}\int_{B(0,1)}{|z|^2 dg}\right),$$
It is easy to verify that this map is admissible.
Since $f$ is assumed to be radial, it follows immediately that
$$ \int_{B(0,1)}{m(x)|f(x)|^2dg} = (\textbf(0)_{\mathbb{R}^d}, t)$$
where the real variable $t$ is easily seen to satisfy
$$-\frac{1}{|B(0,1)|}\int_{B(0,1)}{|z|^2 dg} \leq t \leq 1- \frac{1}{|B(0,1)|}\int_{B(0,1)}{|z|^2 dg}.$$
\begin{center}
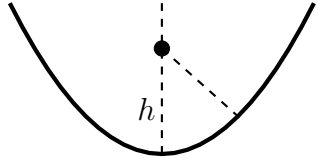
\begin{figure}[h!]
\begin{tikzpicture}[scale = 2]
\draw[ultra thick, domain=-1:1] plot (\x, {\x*\x});
\draw[thick, dashed] (0,0)--(0,1);
\draw [fill] (0,0.7) circle [radius=0.05];
\draw[thick, dashed] (0,0.7)--(0.5,0.5^2);
\draw (-0.1,0.3) node{\Large $h$};
\end{tikzpicture}
\caption{For points of the type $(\textbf(0)_{\mathbb{R}^d}, t)$ the shortest distance to the parabola is comparable to the height $h$.}
\end{figure}
\end{center}
A simple calculation in one dimension shows that for all $0 \leq h \leq 1$
$$ \inf_{-1 \leq x \leq 1}{\|(0,h) - (x,x^2)\|_{\ell^2(\mathbb{R}^2)}} \geq  \frac{\sqrt{3}}{2}\|(0,h)\|_{\ell^2(\mathbb{R}^2)} =  \frac{\sqrt{3}}{2}h$$
and thus, we can deduce that
\begin{align*} \inf_{z \in M}{\left\|m(z) - \int_{M}{m(x)|f(x)|^2dg}\right\|_{\ell^2(\mathbb{R}^{d+1})}}
&\geq \frac{\sqrt{3}}{2}\left\|m(\textbf{0}) - \int_{M}{m(x)|f(x)|^2dg}\right\|_{\ell^2(\mathbb{R}^{d+1})} \\
&= \frac{\sqrt{3}}{2}\int_{B(0,1)}{|x|^2|f(x)|^2dx}.
\end{align*}

Since $f$ has 99\% of its $L^2-$mass at distance at most $1/10$ from the origin, a simple calculation shows that the third term is bounded uniformly away from 0 and we are left with $$\int_{B(0,1)}{|x|^2|f(x)|^2dx}\int_{B(0,1)}{|\nabla f(x)|^2dx} \gtrsim_d 1.$$
\end{proof}

\subsection{The more general inequality}The proof in the last section has a curious structure: we used the symmetry of the function $f$ to bound one of the three terms in our main inequality from below. For general nonradial $f$, the term can be much smaller than the rough estimate we used. Indeed, this
means that for general $f$ our inequality becomes stronger than the classical Euclidean uncertainty principle centered at the origin. We will now understand how that happens and
 show that our result actually behaves like the inequality
$$\inf_{a \in \mathbb{R}^d}\left(\int_{\mathbb{R}^d}{|x-a|^2 f^2 dx}\right)\left(\int_{\mathbb{R}^d}{|\nabla f|^2dx}\right) \gtrsim_d \|f\|_{L^2(\mathbb{R}^d)}^4,$$
which satisfies Breitenberger's suggested property of invariance.
Suppose we are dealing with a function exhibiting spatial
concentration around some point $\textbf{0} \neq a \in B(0,1)$; for simplicity, assume that $f_{\varepsilon}(x)^2$ is given by the probability density of a Gaussian centered around $a$ with variance $\varepsilon$. We will study all three terms in our main inequality
seperately; the gradient term 
$$ \int_{B(0,1)}{|\nabla f_{\varepsilon}|^2dg} \quad \mbox{scales precisely as it would in the Euclidean space for $\varepsilon$ small.}$$
The second term is now of great interest
$$ \inf_{z \in B(0,1)}{\left\|m(z) - \int_{B(0,1)}{m(x)|f_{\varepsilon}(x)|^2dg}\right\|_{\ell^2(\mathbb{R}^d)}} \sim \left\|(a,|a|^2)- \int_{B(0,1)}{(x,|x|^2)|f_{\varepsilon}(x)|^2dg}\right\|_{\ell^2(\mathbb{R}^d)}.$$
We use a change of coordinate to move everything to the origin, expand the weight and use the symmetry of the function to get a cancellation of the linear term 
\begin{align*} \int_{B(0,1)}{(x,|x|^2)|f_{\varepsilon}(x)|^2dg} &\sim \int_{B(0,1-|a|)}{(a+x,|a|^2+2\left\langle a, x \right\rangle + |x|^2)|f_{\varepsilon}(a+x)|^2dg}\\
&\sim (a,|a|^2) + \int_{B(0,1-|a|)}{(x,|x|^2)|f_{\varepsilon}(a+x)|^2dg}
\end{align*} 
The constant vector gets subtracted and what remains is precisely the classical Euclidean uncertainty principle
centered at $a$ since, for $\varepsilon$ small,
\begin{align*}
\left\| \int_{B(0,1-|a|)}{(x,|x|^2)|f_{\varepsilon}(a+x)|^2dg}\right\|_{\ell^2(\mathbb{R}^d)} &\sim  \left\| \left(\textbf{0}_{\mathbb{R}^d}, \int_{B(0,1-|a|)}{|x|^2|f_{\varepsilon}(a+x)|^2dg}\right)\right\|_{\ell^2(\mathbb{R}^d)}\\
&=  \int_{B(0,1)}{|x-a|^2|f_{\varepsilon}(x)|^2dg}
\end{align*}
The third term will behave like a constant as $\varepsilon$ becomes small (the precise value depends on the center of gravity of the segment of the paraboloid and its distance to $a$).

\section{Proof of the Theorem}
\subsection{Outline.} We aim to show that the product of three quantities can never be very small. We start by nothing that for very simple reasons one of the three terms can actually never be too small: the Lipschitz continuity
of $m$ combined with the third condition implies that
$$ \sup_{z \in M}{\|m(z)\|_{\ell^2}} \lesssim_{} \mbox{diam}(M)L \lesssim_{(M,g)} L$$
and therefore
$$\left\|\int_{M}{m(x)|f(x)|^2dg(x)} \right\|_{\ell^2(\mathbb{R}^{d})}^{-2} \gtrsim_{(M,g)} \frac{1}{L^2}.$$
This is not altogether surprising since the purpose of that quantity is to be very large for functions being
very close to constant (and thus yielding a small gradient). It remains to study the two remaining cases,
where either 
$$  \int_{M}{|\nabla f|^2dg} \quad \mbox{or} \quad \inf_{z \in M}{\left\|m(z) - \int_{M}{m(x)|f(x)|^2dg}\right\|_{\ell^2(\mathbb{R}^d)}} \quad \mbox{is small.}$$

\subsection{Case 1 (Frequency Concentration)} An elementary formulation of a Poincar\'{e}-type inequality on connected manifolds is that there exists some $\varepsilon_1 > 0$ depending only on $(M,g)$, such that any $L^2-$normalized function $f$ satisfying 
$$  \int_{M}{|\nabla f|^2dg} \leq \varepsilon_1 \quad \mbox{satisfies the uniform lower bound} \quad
f(x) \geq \frac{1}{2|M|^{\frac{1}{2}}}.$$
We assume now that
$$ \int_{M}{|\nabla f|^2dg} = \varepsilon \leq\varepsilon_1 ~ \mbox{and show that}~~ 
\left\|\int_{M}{m(x)|f(x)|^2dg} \right\|_{\ell^2(\mathbb{R}^{d})} ~~ \mbox{has to be large.}$$
We subtract the mean value of $f$ and consider
$$h = f - \frac{1}{|M|}\int_{M}{f(x) dg} \quad \mbox{satisfying} \quad \int_{M}{|\nabla h|^2dg} = \varepsilon.$$
Since $h$ has mean value 0, we may thus use the Poincar\'{e} inequality for functions satisfying a Neumann condition to conclude
$$ \int_{M}{h^2 dg} \leq g_1^{-1}\int_{M}{|\nabla h|^2dg} = \mu_1^{-1}\varepsilon \lesssim_{(M,g)} \varepsilon,$$
where $\mu_1 > 0$ is the first nontrivial eigenvalue of the Neumann-Laplacian on $M$. We substitute
$$\int_{M}{m(x)f(x)^2dg} = \int_{M}{m(x)\left(\frac{1}{|M|}\int_{M}{f(x) dg} + h(x)\right)^2dg}$$
and, after expanding the square, note that the first term vanishes completely due to the normalization condition
on the map $m:M\rightarrow \mathbb{R}^d$
$$  \int_{M}{m(x)\left(\frac{1}{|M|}\int_{M}{f(x) dg}\right)^2 dg} =
\left(\frac{1}{|M|}\int_{M}{f(x) dg}\right)^2 \int_{M}{m(x)dg} = \textbf{0} \in \mathbb{R}^d.$$
The third term in the expansion also has to be small, pulling out the $L^{\infty}-$norm of $m$ yields
$$ \left\|\int_{M}{m(x) h(x)^2 dg}\right\|_{\ell^2(\mathbb{R}^d)} \leq \left\|\left\|m\right\|_{\ell^2(\mathbb{R}^d)}\right\|_{L^{\infty}(M)}\int_{M}{h(x)^2dg} \lesssim \mbox{diam}(M) L \mu_1^{-1} \varepsilon \lesssim_{(M,g)} L \varepsilon.$$
The real contribution to the size of the entire term comes from the mixed term; using H\"older's inequality,
we first show that the absolute constant is small
$$ \left| \frac{1}{|M|}\int_{M}{f(x) dg}\right| \leq \frac{1}{|M|}\left(\int_{M}{f(x)^2 dg}\right)^{\frac{1}{2}}|M|^{\frac{1}{2}} \leq \frac{1}{|M|^{\frac{1}{2}}} \lesssim_M 1.$$
It remains to estimate the linear term, using the previous estimate and H\"older's inequality once again,
we end up getting 
\begin{align*}\left|\frac{2}{|M|}\int_{M}{f(x) dg}\right|\left\| \int_{M}{m(x)h(x)dg}\right\|_{\ell^2(\mathbb{R}^d)}
&\lesssim  \left\|\left\|m\right\|_{\ell^2(\mathbb{R}^d)}\right\|_{L^{\infty}(M)}\left(\int_{M}{h(x)^2dg}\right)^{\frac{1}{2}} \\
&\lesssim \mbox{diam}(M) L \sqrt{\varepsilon}.
\end{align*}
Collecting all these estimates, we see that if
$$ \int_{M}{|\nabla f|^2dg} = \varepsilon \qquad \mbox{then} \qquad 
\left\| \int_{M}{m(x)f(x)^2dg}\right\|_{\ell^2(\mathbb{R}^d)} \lesssim_{(M,g)} \mbox{diam}(M) L (\sqrt{\varepsilon} + \varepsilon)$$
and thus
$$ \left(\int_{M}{|\nabla f|^2dg}\right)\left\|\int_{M}{m(x)|f(x)|^2dg(x)} \right\|_{\ell^2(\mathbb{R}^{d})}^{-2} \gtrsim_{(M,g)} \frac{1}{L^2 \mbox{diam}(M)^2} \gtrsim_{(M,g)} \frac{1}{L^2}.$$
This shows that the product of two out of the three terms cannot possibly be too small, however, 
one term is left. We need to show that if we have a small gradient, then the center of mass cannot
be close to any element of $m(M) \subset \mathbb{R}^d$. We will now prove that under the assumption above
$$ \inf_{z \in M}{\left\|m(z) - \int_{M}{m(x)|f(x)|^2dg}\right\|_{\ell^2(\mathbb{R}^d)}} \gtrsim_{(M,g)} \frac{C}{L^2}.$$
From the geometric condition on the map $m$, we can conclude that for all $z \in M$
$$\left\|m(z) - \int_{M}{m(x)|f(x)|^2dg}\right\|_{\ell^2(\mathbb{R}^d)} \geq C\int_{M}{ f(x)^2\|m(z) -m(x)\|_{L^2(\mathbb{R}^d)}^2}d g.$$
Using the bilipschitz property of $m$, we get
$$ C\int_{M}{ f(x)^2\|m(z) -m(x)\|_{L^2(\mathbb{R}^d)}^2  dg} \geq \frac{C}{L^2}\int_{M}{f(x)^2 d(z,x)^2dg}.$$
However, since $\varepsilon \leq \varepsilon_1$ and the function is almost constant, we have that
$$ \int_{M}{f(x)^2 d(z,x)dg} \geq \frac{1}{4|M|}\int_{M}{d(z,x)^2dg}  \gtrsim_{(M,g)} 1.$$

\subsection{Case 2 (Spatial Concentration)} It follows from the definition of the manifold that it
behaves essentially like Euclidean space on small scales; in particular, for any compact manifold
there exists some $\varepsilon_2 > 0$ such that for any $z \in M$ and any function containing 
at least 99\% of its $L^2-$mass within a ball of radius $r \leq \varepsilon_2$, then we have the
classical uncertainty principle
$$\int_{M}{f(x)^2 d(z,x)^2dg(x)} \int_{M}{|\nabla f(x)|^2 dg} \gtrsim_{(M,g)} 1.$$
Assume now that for some fixed $z \in M$
$$ \varepsilon := \left\|m(z) - \int_{M}{m(x)|f(x)|^2dg}\right\|_{\ell^2(\mathbb{R}^d)} \leq \frac{\varepsilon_2^2 C}{100 L^2}.$$
We will now be using the curvature condition of $m$ as follows: it states that
for some $0 < C < \infty$, for all
integers $N \geq 1$ and all elements $x_1, \dots, x_N,z \in M$:
$$ \left\| \frac{1}{N}\sum_{i=1}^{N}{m(x_i)} - m(z)\right\|_{\ell^2(\mathbb{R}^d)} \geq 
 \frac{C}{N}\sum_{i=1}^{N}\left\|m(x_i) - m(z) \right\|^2_{\ell^2(\mathbb{R}^d)}.$$
If we now let $N \rightarrow \infty$ and position points in such a way that their empirical density converges weakly to
$f(x)^2$, then we have
\begin{align*}
\lim_{N \rightarrow \infty}{ \left\| \frac{1}{N}\sum_{i=1}^{N}{m(x_i)} - m(z)\right\|_{\ell^2(\mathbb{R}^d)}}&= \left\| \int_{M}{f(x)^2m(x)dg} - m(z)\right\|_{\ell^2(\mathbb{R}^d)} \\
\lim_{N \rightarrow \infty}{ \frac{C}{N}\sum_{i=1}^{N}\left\|m(x_i) - m(z) \right\|^2_{\ell^2(\mathbb{R}^d)}}&= C\int_{M}{f(x)^2 d(m(z),m(x))^2dg}
\end{align*}

Using again the geometric condition satisfied by $m$ and the bilipschitz property, we conclude
\begin{align*}
 \varepsilon =\left\|m(z) - \int_{M}{m(x)|f(x)|^2dg}\right\|_{\ell^2(\mathbb{R}^d)}
&\geq C\int_{M}{f(x)^2 d(m(z),m(x))^2dg} \\
&\geq\frac{C}{L^2}\int_{M}{f(x)^2 d(z,x)^2dg}
\end{align*}
Then Markov's inequality immediately implies that
$$ \int_{B(z,\varepsilon_2)}{f(x)^2 dx} \geq \frac{99}{100}.$$
Then, however,
$$\int_{M}{f(x)^2 d(z,x)^2dg} \int_{M}{|\nabla f(x)|^2 dg} \gtrsim_{(M,g)} 1.$$
As before, we may use estimate the remaining third term trivially with
$$\left\|\int_{M}{m(x)|f(x)|^2dg} \right\|_{\ell^2(\mathbb{R}^{d})}^{-2} \geq
\frac{1}{\mbox{diam}(M)^2 L^2} \gtrsim_{(M,g)} \frac{1}{L^2}.$$

\subsection{Conclusion of the argument.} Suppose neither Case 1 nor Case 2 applies. Then
$$  \int_{M}{|\nabla f|^2dg} \geq \varepsilon_1 \quad \mbox{and} \quad
\inf_{z \in \mathbb{T}}{\left\|m(z) - \int_{\mathbb{T}}{m(x)|f(x)|^2dg}\right\|_{\ell^2(\mathbb{R}^2)}} \geq  \frac{\varepsilon_2^2 C}{100 L^2},$$
where $\varepsilon_1, \varepsilon_2$ are universal constants depending only on the manifold. Using again
$$\left\|\int_{M}{m(x)|f(x)|^2dg(x)} \right\|_{\ell^2(\mathbb{R}^{d})}^{-2} \geq
\frac{1}{\mbox{diam}(M)^2 L^2} \gtrsim_{(M,g)} \frac{1}{L^2},$$
the result follows. $\qed$

\section{Disconnected manifolds}
The problem is of a very different nature once the compact manifold is disconnected. The proof
can essentially be adapted but we will lose control over the behavior of the implicit constant. It is quite interesting that the degeneracy induced by disconnected manifolds can be observed in our main statement as follows: take a connected manifold and let a small part of it degenerate. The constant in our main theorem did depend (among many other things) on the inverse of the first eigenvalue
of the Dirichlet-Laplacian $g_1$: this quantity tends to 0 in the process of neck-pinching and $g_1(M)=0$ in the case of a disconnected manifold.
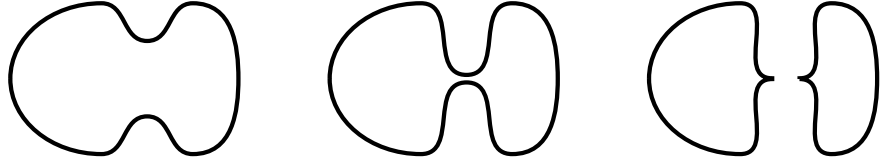
\begin{figure}[h!]
\begin{center}
\begin{tikzpicture}[xscale = 0.6,yscale = 0.5]
\draw [ultra thick] (0,0) to [out=270,in=180] (2,-2) to [out=0,in=180] (3,-1) to [out=0,in=180] (4,-2) to [out=0, in=270] (5,0) to [out=90, in = 0] (4,2) to [out=180, in = 0] (4,2) to [out=180, in = 0] (3,1) to [out=180, in = 0] (2,2) to [out=180, in = 90] (0,0);
\draw [ultra thick] (7,0) to [out=270,in=180] (9,-2) to [out=0,in=180] (10,-0.1) to [out=0,in=180] (11,-2) to [out=0, in=270] (12,0) to [out=90, in = 0] (11,2) to [out=180, in = 0] (11,2) to [out=180, in = 0] (10,0.1) to [out=180, in = 0] (9,2) to [out=180, in = 90] (7,0);
\draw [ultra thick] (14,0) to [out=270,in=180] (16,-2) to [out=0,in=180] (16.7,-0.01) to [out=90,in=270]  (16.7,0.01) to [out=180, in = 0] (16,2) to [out=180, in = 90] (14,0);
\draw [ultra thick](17.3,-0.01) to [out=0,in=180] (18,-2) to [out=0, in=270] (19,0) to [out=90, in = 0] (18,2) to [out=180, in = 0] (18,2) to [out=180, in = 0] (17.3,0.01) to [out=270, in = 90] (17.3,-0.01);
\end{tikzpicture}
\end{center}
\caption{Neck pinching of a manifold.}
\end{figure}

\subsection{Setup}  Let $(M,g)$ be some compact $n-$dimensional, smooth manifold (without boundary or with smooth boundary)
and $k$ connected components
$$ M = M_1 \cup M_2 \cup \dots \cup M_k.$$
 We demand that the map $m: M \rightarrow \mathbb{R}^d$ be continuous and bijective and that it satisfies a uniform convexity condition with $0 < C < \infty$
$$ \forall~N \in \mathbb{N}, x_1, \dots, x_N,z \in M \quad \left\| \frac{1}{N}\sum_{i=1}^{N}{m(x_i)} - m(z)\right\|_{\ell^2(\mathbb{R}^d)} \geq 
 \frac{C}{N}\sum_{i=1}^{N}\left\|m(x_i) - m(z) \right\|^2_{\ell^2(\mathbb{R}^d)}.$$
We define $k$ points $p_1, \dots, p_k \in \mathbb{R}^d$ via
$$ p_i = \int_{M_i}{m(x)dg}$$
and the simplex $S \subset \mathbb{R}^d$ via
$$ S=\left\{\sum_{i=1}^{k}{\alpha_i p_i}: \alpha_i \geq 0 \wedge \sum_{i=1}^{k}{\alpha_i}=1\right\}.$$

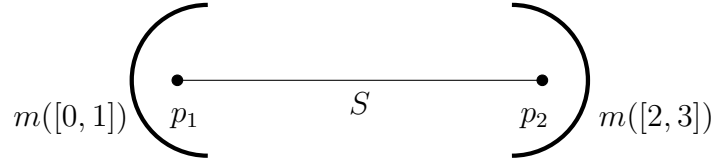
\begin{figure}[h!]
\begin{center}
\begin{tikzpicture}
\draw [ultra thick] (0,0) to [out=180,in=90] (-1,-1) to [out=270,in=180] (0,-2);
\draw (-1.8,-1.5) node{\Large $m([0,1])$};
\draw [ultra thick] (4,0) to [out=0,in=90] (5,-1) to [out=270,in=0] (4,-2);
\draw (5.9,-1.5) node{\Large $m([2,3])$};
\draw [fill] (-0.4,-1) circle [radius=0.07];
\draw (-0.3,-1.5) node{\Large $p_1$};
\draw [fill] (4.4,-1) circle [radius=0.07];
\draw (4.3,-1.5) node{\Large $p_2$};
\draw (2,-1.3) node{\Large $S$};
\draw (4.4,-1) -- (-0.4, -1);
\end{tikzpicture}
\end{center}
\caption{The imbedding of the manifold $M = [0,1] \cup [2,3]$ into $\mathbb{R}^2$ via a suitable map
$m$, the points $p_1, p_2$ and the associated simplex $S$ (here a line).}
\end{figure}
We assume furthermore that $m$ is bilipschitz with constant $L$ on every connected component $M_i$. In order for our main argument to work, we require additionally some control on the separation of the components of $M$ in the mapping $m$. We define the separation distance $\sigma$ via
$$ \sigma = \min_{1 \leq i < j \leq k}{\quad \min_{x \in M_i, y \in M_j}{\quad \|m(x) - m(y)\|_{\ell^2(\mathbb{R}^d)}}}.$$ Since all $M_i$ are compact and $m$ is bijective as well as continuous, we always have $\sigma > 0$. However, the implicit constant in our statement will also depend on $\sigma$ -- it is not clear to us whether this is really necessary or whether the statement could be stated without any control on $\sigma$.
\begin{thm} Let $(M,g)$ and $m$ be as above. For any real-valued $f \in H^1(M)$ with $\|f\|_{L^2(M)} = 1$ we define the quantity
$$ U = \left(\inf_{z \in M}{\left\|m(z) - \int_{M}{m(x)|f(x)|^2dg}\right\|_{\ell^2(\mathbb{R}^d)}}\right)
\left(\inf_{s \in S}{\left\|s - \int_{M}{m(x)|f(x)|^2dg}\right\|_{\ell^2(\mathbb{R}^d)}^{-2}}\right) \int_{M}{|\nabla f|^2dg}.$$
Then
$$ U \gtrsim_{(M,g), C, L, \sigma}1.$$
\end{thm}
The proof is a rather straightforward adaption of our main argument; we will only 
sketch the necessary modifications and emphasize the arising difficulty where the argument starts to lose
track of the constant.

\begin{proof} 
\textbf{Case 1 (Frequency Concentration).} 
The case of frequency localization is only slightly different: if the function has a small gradient, then it is essentially almost constant on every connected component $M_i$.
The way the simplex was defined guarantees that the weighted average of the respective center of masses
is close to the simplex and that the center of mass is therefore close to the simplex independent of the actual numerical value of the function on each connected component. In greater detail, we may say that if
$$ \int_{M}{|\nabla f|^2 dg} = \varepsilon \ll 1 \qquad \mbox{is small,}$$
then it is trivially also small on all the pieces $M_i$. We may thus run our previous argument on each $M_i$ and may deduce, in the same way as before, that then
$$ \left\| \int_{M_i}{m(x)|f(x)|^2 dg} - p_i \int_{M_i}{f(x)^2 dg}\right\|_{\ell^2(\mathbb{R}^d)} \lesssim_{(M,g)} L\sqrt{\varepsilon},$$
where the implicit constant depends, among other things, on the first eigenvalue of the Neumann-Laplacian on $M_i$. Note that
$$ \sum_{i=1}^{k}{p_i \int_{M_i}{f(x)^2 dg}} \in S$$
and therefore, using the triangle inequality,
\begin{align*} \inf_{s \in S}{\left\|s - \int_{M}{m(x)|f(x)|^2dg}\right\|_{\ell^2(\mathbb{R}^d)}^{-2}}& \geq
\left\|\sum_{i=1}^{k}{p_i \int_{M_i}{f(x)^2 dg}}  - \int_{M}{m(x)|f(x)|^2dg}\right\|_{\ell^2(\mathbb{R}^d)}^{-2}\\
&\gtrsim_{(M,g)} \frac{1}{L^2 \varepsilon}.
\end{align*}
The bound
$$ \inf_{z \in M}{\left\|m(z) - \int_{M}{m(x)|f(x)|^2dg}\right\|_{\ell^2(\mathbb{R}^d)}} \gtrsim_{(M,g)} \frac{C}{L^2}$$
can be obtained as before without any modifications. \\

\textbf{Case 2 (Mass Concentration).}
The case of mass concentration is somewhat different: let us suppose that again that
$$\left\|m(z) - \int_{M}{m(x)|f(x)|^2dg}\right\|_{\ell^2(\mathbb{R}^d)} = \varepsilon \qquad \mbox{is small}.$$
Before we used the geometric condition satisfied by $m$ and the bilipschitz property to conclude that
\begin{align*}
 \varepsilon =\left\|m(z) - \int_{M}{m(x)|f(x)|^2dg}\right\|_{\ell^2(\mathbb{R}^d)} 
&\geq C\int_{M}{f(x)^2\|m(z) - m(x)\|_{L^2(\mathbb{R}^d)}^2}dg\\
&\geq \frac{C}{L^2}\int_{M}{f(x)^2 d(z,x)^2dg}.
\end{align*}
The last step in the inequality will now fail because we only have the bilipschitz property of $m$ on every connected component of the manifold but not globally on $M$. However, the geometric condition on $m$ remains valid and allows to conclude that at least 
$$ \varepsilon = \left\|m(z) - \int_{M}{m(x)|f(x)|^2dg}\right\|_{\ell^2(\mathbb{R}^d)} 
\geq C\int_{M}{\|m(z) - f(x)^2m(x)\|_{\ell^2(\mathbb{R}^d)}^2}dg.$$
If there was some way to conclude from this that at least $99\%$ of the $L^2-$mass of $f$ is contained within the one connected component containing $z$, we could again 
emulate the same argument as before. This is where the control on the seperation of the images $m(M_i)$ enters; if it were not the case that at least  $99\%$ of the $L^2-$mass 
of $f$ is contained within one connected component containing $m(z)$, then we could conclude that
$$ C\int_{M}{\|m(z) - f(x)^2m(x)\|_{L^2(\mathbb{R}^d)}^2}dg \geq \frac{C \sigma^2}{100},$$
which is a contradiction for $\varepsilon$ sufficiently small (depending only on $\sigma$). 
\end{proof}

\section{Open problems and remarks}

\subsection{An inverse problem.} In this section we will be concerned with the question how one would go about picking the right map $m$ to obtain 
as much information about $f$ as possible. Doing the same sort of simplification we did before, our problem may be stated as follows.
\begin{quote} \textbf{Inverse problem.} Given $(M,g)$ and $f:M \rightarrow \mathbb{R}$, which map $m:M \rightarrow \mathbb{R}^d$ is both very curved ($C$ big) as well as
isometric as possible ($L$ small) but also has the property that when it is weighted with $|f(x)|^2$ yields a point far away from the origin? In short,
prove lower bounds
$$\sup_{m:M \rightarrow \mathbb{R}^d}\frac{C}{L^5}\left\| \int_{M}{m(x)|f(x)|^2dg}\right\|_{\ell^2(\mathbb{R}^d)}^{2} \qquad \mbox{in terms of}~f$$
and describe the mappings $m$ achieving that bound.
\end{quote}
Recall that the motivation for studying this term comes from our main result, which implies that
$$  \int_{M}{|\nabla f|^2dg} \gtrsim_{(M,g),d} \frac{C}{L^5}\left\| \int_{M}{m(x)|f(x)|^2dg}\right\|_{\ell^2(\mathbb{R}^d)}^{2}.$$

We will give a sample result for the case $m:\mathbb{T} \rightarrow \mathbb{R}^3$, where the torus has length $2\pi$ -- gaining some insight into how such a result
could be sharpened seems to be an interesting problem.
\begin{proposition} Let $f \in C^{1}(\mathbb{T})$ be normalized in $L^2(\mathbb{T})$. We have
$$\sup_{m:\mathbb{T} \rightarrow \mathbb{R}^3}\frac{C}{L^5}\left\|\int_{\mathbb{T}}{m(x)|f(x)|^2dg} \right\|^2_{\ell^2(\mathbb{R}^{3})} \gtrsim \frac{1}{(1+\|\nabla f\|_{L^{\infty}})^7}\left( \int_{\mathbb{T}}{|f(x)|^4dx}
- \frac{1}{2\pi}\right) \geq 0.$$
\end{proposition}
H\"older's inequality immediately implies that the last factor only vanishes for constant functions: we will
therefore always recover some nontrivial information.

\begin{proof} The proof is constructive: consider the map $m:\mathbb{T} \rightarrow \mathbb{R}^3$ with
$$ m(t) = (\cos{t}, \sin{t}, |f(t)|^2 - \frac{1}{2\pi}).$$
We start by estimating the Lipschitz constant; clearly, for $s,t \in \mathbb{T}$
$$|s-t| \lesssim \|m(s) - m(t)\|_{\ell^2(\mathbb{R}^3)} \leq (1+\|\nabla f\|_{L^{\infty}})|s-t|.$$
This implies that $L \lesssim  (1+\|\nabla f\|_{L^{\infty}})$.
To estimate $C$ (which will depend strongly on $f$), we use a simple trick. Recall that we have
to find a constant $C$ such that the inequality
$$ \left\| \frac{1}{N}\sum_{i=1}^{N}{m(x_i)} - m(z)\right\|_{\ell^2(\mathbb{R}^3)} \geq 
 \frac{C}{N}\sum_{i=1}^{N}\left\|m(x_i) - m(z) \right\|^2_{\ell^2(\mathbb{R}^3)}$$
holds. We use $\pi:\mathbb{R}^3 \rightarrow \mathbb{R}^2$ to denote the projection
onto the first 2 coordinates, i.e. $\pi(x,y,z) = (x,y)$. Trivially,
\begin{align*} \left\| \frac{1}{N}\sum_{i=1}^{N}{m(x_i)} - m(z)\right\|_{\ell^2(\mathbb{R}^3)}
&\geq  \left\| \pi\left(\frac{1}{N}\sum_{i=1}^{N}{m(x_i)}\right) - \pi(m(z))\right\|_{\ell^2(\mathbb{R}^2)} \\
&= \left\| \frac{1}{N}\sum_{i=1}^{N}{\pi\left(m(x_i)\right) }- \pi(m(z))\right\|_{\ell^2(\mathbb{R}^2)}
\end{align*}
The canonical embedding of $\mathbb{T}$ into $\mathbb{R}^2$, i.e. $x \rightarrow (cos(x), \sin(x))$ satisfies
the geometric curvature condition for some constant and therefore
$$ \left\| \frac{1}{N}\sum_{i=1}^{N}{\pi\left(m(x_i)\right) }- \pi(m(z))\right\|_{\ell^2(\mathbb{R}^2)}
\gtrsim \frac{1}{N}\sum_{i=1}^{N}{ \left\| \pi\left(m(x_i)\right) - \pi(m(z))\right\|^2_{\ell^2(\mathbb{R}^2)}}.$$
We conclude the argument by using the already derived bound on the Lipschitz constant to deduce that
$$ (1+L)^2 \left\| \pi\left(m(x_i)\right) - \pi(m(z))\right\|^2_{\ell^2(\mathbb{R}^2)} \gtrsim 
\left\| m(x_i) - m(z)\right\|^2_{\ell^2(\mathbb{R}^3)}$$
Altogether, this implies that
\begin{align*} \left\| \frac{1}{N}\sum_{i=1}^{N}{m(x_i)} - m(z)\right\|_{\ell^2(\mathbb{R}^3)}
&\gtrsim \frac{1}{N}\sum_{i=1}^{N}{ \left\| \pi\left(m(x_i)\right) - \pi(m(z))\right\|^2_{\ell^2(\mathbb{R}^2)}} \\
&\gtrsim \frac{1}{(1+L)^2}\frac{1}{N}\sum_{i=1}^{N}{ \left\| (m(x_i) - m(z)\right\|^2_{\ell^2(\mathbb{R}^3)}} 
\end{align*}
and this implies that
$$ C \gtrsim \frac{1}{ (1+\|\nabla f\|_{L^{\infty}})^2}.$$
For the conclusion of the argument, we use $\pi_3:\mathbb{R}^3 \rightarrow \mathbb{R}$ to denote
the projection onto the third component $\pi_3(x,y,z) = z$. The Phytagorean theorem and H\"older's inequality yield that
$$\left\|\int_{\mathbb{T}}{m(x)|f(x)|^2dg} \right\|_{\ell^2(\mathbb{R}^{2})}^2 \geq \left|\pi_3\left(\int_{\mathbb{T}}{m(x)|f(x)|^2dg}\right)\right| = \int_{\mathbb{T}}{|f(x)|^4dx} - \frac{1}{2\pi} \geq 0$$
with equality if and only if $f$ is constant. This gives the statement.
\end{proof}

\subsection{A weaker geometric condition}  Simple examples show that the curvature condition on $m$ is necessary. However, certainly some statements remain true if the geometric condition is allowed to be violated
in a localized manner in a controlled way. We give a very simple example in the setting $m:[-1,1] \rightarrow \mathbb{R}^2$ with
$$m(x) = \left(x,|x|^k - \frac{1}{k+1}\right)\quad \mbox{for some integer} ~k \geq 3.$$ 
The map satisfies the geometric condition outside of any neighbourhood of the origin; however, the curvature
decays as $\kappa \sim |x|^{k-2}$ as we approach the origin.
\begin{proposition} We have
$$ \left(\inf_{-1 \leq z \leq 1}{\left\|m(z) - \int_{-1}^1{m(x)|f(x)|^2dx}\right\|^{\frac{2}{k}}_{\ell^2(\mathbb{R}^2)}}\right) \left\| \int_{-1}^{1}{m(x)|f(x)|^2dx}\right\|_{\ell^2(\mathbb{R}^2)}^{-2} \int_{-1}^1{|\nabla f|^2dx} \gtrsim_k 1.$$
\end{proposition}
The proof follows immediately from repeating our argument in the non-degenerate case and adjusting
for different parameters in the case of spatial localization around the origin. Of course, the case
$k=2$ reduces to a special case of our result.

\subsection{A purely combinatorial problem.} This section is about a problem of a more geometric-combinatorial type that arises somewhat naturally in this context. Throughout the paper we dealt with mappings satisfying
 for some $0 < C < \infty$
$$ \forall~N \in \mathbb{N}~\forall~ x_1, \dots, x_N,z \in M \quad \left\| \frac{1}{N}\sum_{i=1}^{N}{m(x_i)} - m(z)\right\|_{\ell^2(\mathbb{R}^d)} \geq 
 \frac{C}{N}\sum_{i=1}^{N}\left\|m(x_i) - m(z) \right\|^2_{\ell^2(\mathbb{R}^d)}.$$
This condition was clearly necessary and it is easy to construct example where its failure implies
the failure of the desired inequality, however, one might ask whether this condition may be concluded
from a simpler (and easier-to-check) condition. It is easy to see that if $M=[0,1]$ or $M=\mathbb{T}$
and $m:M \rightarrow \mathbb{R}^2$, then it suffices to check the condition for $N=2$ because that
already allows to bound the curvature from below, which then implies all necessary results. It is not
difficult to see that $N=2$ is in general not sufficient for maps 
$m:\mathbb{T} \rightarrow \mathbb{R}^3$. In contrast, if $\mathbb{D} = B(0,1) \subset \mathbb{R}^2$ is the unit disk in the plane, then for maps $m:\mathbb{D} \rightarrow \mathbb{R}^3$ it suffices again to
take $N=2$. It seems natural to conjecture that for each fixed setting it suffices to check a certain finite number
of points (these examples suggest that in the case of maps $m:M \rightarrow \mathbb{R}^d$ it might perhaps suffice to check $N = d - $dim($M$)$ + 1$ points). 

\subsection{Limits of the uncertainty principle.} It is easily seen that no mapping $m: \mathbb{T}^2 \rightarrow \mathbb{R}^3$ satisfying the geometric conditions can exist (because for any embedding there will always be a line intersecting $m(\mathbb{T}^2)$ in at least three points. This is not due to a failure of our methods but simply because inequalities of the type we studied do not hold in this case. Are there any natural analogues?\\

\textbf{Acknowledgement.} The author was supported by SFB 1060 of the DFG.

\end{document}